# Soundness does not come for free (if at all)[*]

Kaave Lajevardi

La Société des Philosophes Chômeurs, P.O. Box 19888-3844, Téhéran, IRAN

`kaave.lajevardi@gmail.com`

Saeed Salehi

University of Tabriz, P.O. Box 51666-16471, Tabriz, IRAN

`root@SaeedSalehi.ir`

**Abstract.**

We respond to some of the points made by Bennet and Blanck (2022) in this JOURNAL concerning a previous publication of ours elsewhere (2021).

———

---

[*] *Added 19 October 2023.* As indicated in the Abstract, this is a reply to Bennet and Blanck's reply to a paper of ours. While we published our paper in *Philosophia Mathematica,* their reply appeared in *Theoria* (the referent of 'this JOURNAL' in the Abstract). The referees of *Theoria* voted against our reply, but the editor-in-chief gave us an opportunity of "respond to the reviewer(s)' comments and resubmit a revised version".

As both of us are now busy with other projects (hopefully more constructive than writing rejoinders), and as we find the current piece adequate *as a reply,* we decided not to do any further work on this note and make it available here with no changes.

In a recent paper [henceforth 'BB'] in *Theoria*, Christian Bennet and Rasmus Blanck have discussed parts of a 2021 paper of ours [henceforth 'LS']. We appreciate their demonstration that part of one of our technical results can be derived from more general observations concerning fixed points of arithmetical formulae, though we find their subsequent criticism unconvincing.

## 1. Soundness is not essential to Gödel's theorems.

Recall that, by definition, a theory in the language of arithmetic is *sound* iff all its theorems hold in the standard structure of natural numbers. In his classic paper (1931), Gödel begins his informal introduction with the assumption that theories in question are sound, but then strengthens his results—both with regard to substance and methodology—by dropping the assumption of soundness in favour of of $\omega$-consistency. Just before embarking on the formal and technical part of the paper, he writes

> The purpose of carrying out the above proof with full precision in what follows is, among other things, to replace the second of the assumptions just mentioned by a purely formal and much weaker one. (Gödel 1931, p. 151.)

The 'second of the assumptions just mentioned' is soundness, while the purely formal and *much weaker* one is $\omega$-consistency. Hence, insofar as Gödel's theorems are concerned, officially it is nowhere assumed that any theory whatsoever is sound.

Bennet and Blanck do not demand that the theories *which are subject to the incompleteness theorems* (e.g., PA, ZF) should be sound; rather, they are talking about soundness of subtheories of these in which it is proven that all of the I-am-not-provable sentences are equivalent—to wit, they require a sound "base". Gödel himself does not seem to be interested in isolating a base theory;[1] but, of course, that might be of historical importance only. Let us explore the issue a bit.

Before doing that, however, let us express our dissatisfaction with one of the claims made by BB.

## 2. ¬Pr$_T$(x) by itself does the job.

Bennet and Blanck claim that their Observation 2 "directly gives Theorem 2 of [LS] as a special case". This is misleading at best. Observation 2 of BB has two parts: (i) if a theory $T$ is unsound then every formula has a false fixed point in $T$; conversely, (ii), if every formula has a false fixed point in $T$ then $T$ in unsound.

While (i) does in fact generalize half of Theorem 2 of LS (to the effect that if $T$ is unsound then ¬Pr$_T$(x) has a false fixed point in $T$), the situation is different with (ii). The second half of Theorem 2 of LS shows that if the *particular* formula ¬Pr$_T$(x) has a false fixed point in $T$ then $T$ is unsound; but (ii) of BB infers the unsoundness of a theory from the evidently stronger assumption that *every* formula has a false fixed point in that theory. Hence (ii) is evidently weaker than Theorem 2 of LS. Unfortunately, this has a deteriorating effect on proof offered by BB—a proof consisting of "Immediate."—of our Theorem 2.

Now, back to the main topic.

---

[1] Thus on page 181 of (1931), when enumerating the properties of a system $P$ used in the proof of Theorem VI (i.e., the first incompleteness theorem, p. 171), he says—in his own terminology—that recursive [*rekursive*] relations should be definable *in P*, where there is no mention of base theories.

## 3. The choice of the base theory: Is PA sound?

The final paragraph of BB tries to justify talking about *the* Gödel sentence of an arithmetical theory, where a key element is using a sound base theory. We are not particularly moved *by* the attempted justification, and there are two reasons for this, the first being that we had already moved in that direction! Thus see Remark 5 of LS for a discussion on going "di-theoretically" and using a sound base theory —see also the Appendix to our (2019).

Secondly, Bennet and Blanck talk about "*the* Gödel sentence of *T over* PA", which implies that, in their view, PA is sound. (There is more than a mere implication in the proof of their Observation 2, where they explicitly write "the sound theory of PA".) One may wonder: what is the justification for the assumption that PA is sound?

We are aware that, in common practice, mathematicians assume that PA is sound. Yet, as we put it elsewhere (2019, footnote 7), this seems to be something like an article of faith rather than an argued-for thesis. We certainly *hope* that PA is sound; but do we *know* that it is sound or even consistent? Is the soundness of PA something like an (meta-)axiom? Is it a legitimate proposition to assume in a debate? We are in doubt.

## 4. More on the impropriety of the *the* talk.

We welcome the opportunity to adduce more evidence against the *the* talk and in defence of our pluralist attitude towards Gödelian sentences.

### 4.1. The diversity of what sentences say.

In LS, our strategy was to find an arithmetical theory *T* and two I-am-unprovable sentences, *A* and *B*, such that *A* is true in the standard model of arithmetic while *B* false. In our view, it is rather odd to call each of *A* and *B* '*the* Gödel sentence of *T*' even though they are *T*-equivalent—rather, we suggest that *T* has Gödelian sentence*s* [in the plural], of which *A* and *B* are two. Each such *T* is of course unsound, and, according to BB, should not be trusted as the theory in which the equivalence of *A* and *B* is proven. In this subsection we play in the territory delimited by BB.

Recall that Bennet and Blanck presume PA to be sound. Hence, they would consider each of its well known subtheories $I\Sigma_n$ to be sound. Let Pr be a standard provability predicate of PA. For each natural number *n,* apply the Diagonal Lemma to get a sentence $\gamma_n$ such that $I\Sigma_1 \vdash \gamma_n \leftrightarrow [\mathrm{Con}(I\Sigma_n) \,\&\, \neg \mathrm{Pr}(\#\gamma_n)]$. Since $I\Sigma_1$ is a subtheory of PA, and since $\mathrm{PA} \vdash \mathrm{Con}(I\Sigma_n)$, we have $\mathrm{PA} \vdash \gamma_n \leftrightarrow \neg \mathrm{Pr}(\#\gamma_n)$, so that each $\gamma_n$ is PA-equivalent to Con(PA), and $\gamma_n$s are pairwise PA-equivalent. We gather that, according to BB, each and every $\gamma_n$ may be called '*the* Gödel sentence of PA [over PA]'.

Now, in the eye of PA, the sentence $\gamma_{95}$ says, inter alia, that $I\Sigma_{95}$ is consistent, while $\gamma_{411}$ says, inter alia, that $I\Sigma_{411}$ is consistent (we ignore their self-referential parts), *and these are different things to say*. If the difference is not already obvious, note that, as is well known, each $I\Sigma_n$ proves $\mathrm{Con}(I\Sigma_m)$ for every $m < n$ but no $I\Sigma_n$ proves its own consistency, by Gödel's second incompleteness theorem (*provided that* it is consistent, which surely is the case according to BB). To label $\gamma_{95}$ and $\gamma_{411}$, together with infinitely many other $\gamma_i$s, collectively as 'the Gödel sentence of PA' is, in our view, dangerously close to subscribing to an implausible hyperextensional doctrine about sentences. It is like thinking of both '13 is a prime' and 'the 1279[th] Mersenne number is a prime' as *the* theorem [in the singular] of PA. There are many PA-theorem*s*, and these are two of them.

## 4.2. Beyond arithmetical theories?

Gödel's incompleteness theorems hold for several theories for which no notion of soundness is generally agreed upon. Set theory provides an outstanding example here. One may be attracted to the view of taking the proper class V as something with respect to which soundness is defined; however, even dropping worries concerning the ontological status of V, such a move begs the question against the view that no single universe captures all the truths of set theory—see Hamkins (2012) and references therein. However, even if one accepts the universe view (which we do not), still it is not clear how to define soundness in set theory. At least with regard to some candidates, one encounters the real possibility of two 'sound' theories contradicting each other (say one of them proves CH, the other refutes it).

Of course, BB is interested in a sound *base*, not a sound theory as such. What they would want for applications of incompleteness to set theory is not that a given *set theory* is sound; rather, they would insist on a sound *subtheory* which is resourceful enough to prove the existence and pairwise equivalence of all Gödel sentences. One may hope for a quite minimal piece of arithmetic: if Q does not work, then perhaps $I\Delta_1$ or $I\Sigma_1$, which should be interpreted in the relevant set theory. But it is not clear that this can be done in a unique way in the case of, e.g., Gödel-Bernays set theory (which is, by the way, the very set theory used by Gödel to prove the consistency of CH in the late 1930s).

Hence our title.